\documentclass[11pt,openany]{article}

\usepackage{amssymb}
\usepackage{amsmath}
\usepackage{amsthm}
\usepackage{amsfonts}
\usepackage{mathrsfs}
\usepackage{makeidx}
\usepackage{multicol}

\input{cyracc.def}

\newtheorem{Theorem}{Theorem}[section]
\newtheorem{Lemma}{Lemma}[section]

\newtheorem{Corollary}[Theorem]{Corollary}
\newtheorem{Remark}[Theorem]{Remark}
\newtheorem{Example}{Example}[section]

\newtheorem{Problem}{Problem}[section]
\numberwithin{equation}{section}

\begin{document}
\title{\Large\bf  Further results on normal families of meromorphic functions\\
concerning shared values}

\author{Xiao-Bin Zhang$^a$\thanks{Correspoding author: E-mail: xbzhang1016@mail.sdu.edu.cn(X.B. Zhang); xujunf@gmail.com(J.F. Xu)  }
\quad and \quad  Jun-Feng Xu$^b$ \\
\small{$^a$College of Science, Civil Aviation
University of China, Tianjin 300300, China}\\
\small{$^b$Department of Mathematics, Wuyi University, Jiangmen, Guangdong
529020, China }}
\date{}
\maketitle

\vspace{3mm}

\begin{abstract}
 In this paper,  we prove two normality criteria for families of some functions
 concerning shared values, the results generalize those given by Hu and Meng.
Some examples are given to show the sharpness of our results.
\end{abstract}


\pagenumbering{arabic}

\section{Introduction and main results}

\quad\,\,\,Let $\mathbb{C}$ denote the complex plane and $f (z)$ be a
non-constant meromorphic  function in $\mathbb{C}$. It is assumed that the
reader is familiar with the standard notation used in the Nevanlinna
value distribution theory such as the characteristic function $T(r,
f)$, the proximity function $m(r, f)$, the counting function $N(r,
f)$ (see, e.g. \cite{H:100,Y:111,Y:110}), and $S(r, f)$ denotes any quantity that satisfies the condition
$S(r, f)=o(T(r, f))$ as $r\to \infty $ outside of a possible
exceptional set of finite linear measure.

\vspace{3mm}
Let $f(z)$ and $g(z)$ be two non-constant meromorphic functions, $a$ be a finite complex number, if $f-a$
and $g-a$ have the same zeros (Ignoring multiplicities), then we say that $f$ and $g$ share $a$.

\vspace{3mm}

 Let $\mathcal {F}$ be a
family of meromorphic functions defined in a domain $D
\subset \mathbb{C}.$ $\mathcal {F}$ is said to be normal in
$D$, in the sense of Montel, if for any sequence $f_{n}
\in \mathcal {F}$, there exists a subsequence $f_{n_{j}}$ such that
$f_{n_{j}}$ converges spherically locally uniformly in $D$, to a meromorphic function or $\infty$ (see \cite{H:100,Y:110}).

\vspace{3mm}
According to Bloch's principle, every condition which reduces a meromorphic
function in $\mathbb{C}$ to a constant, makes a family of meromorphic functions in a domain $D$ normal.
Although the principle is false in general, many authors proved normality criteria for
families of meromorphic functions by starting from Picard type theorems. For instance,
\vspace{3mm}

\noindent{\textbf{Theorem A.}} \cite{H:99}
 \emph{Let $n\geq 5$ be an integer, $a$, $b \in\mathbb{C}$ and $a\neq 0$. If, for a meromorphic function $f$ ,
$f' +af ^{n} \neq b$ for all $z \in\mathbb{C}$, then $f$ must be a constant.
 }

 \vspace{3mm}

\noindent{\textbf{Theorem B.}} \cite{P:108,P:109}
 \emph{Let $n\geq3$ be an integer, $a$, $b \in\mathbb{C}$, $a \neq0$ and $\mathcal {F}$ be a family of meromorphic functions in a domain $D$. If
$f' +af ^{n} \neq b$ for all $f \in\mathcal {F}$, then $\mathcal {F}$ is a normal family.
 }

 \vspace{3mm}

In 2008, Zhang \cite{Z:104} improved theorem B by the idea of shared values, he got

\vspace{3mm}
\noindent{\textbf{Theorem C.}}
 \emph{Let $\mathcal {F}$ be a family of meromorphic functions in D, n be a positive integer and $a$, $b$ be two constants such
that $a\neq 0, \infty$ and $b \neq\infty$. If $n\geq4$ and for each pair of functions $f$ and $g \in \mathcal {F}$, $f'-af^{n}$ and $g'-ag^{n}$ share the
value $b$, then $\mathcal {F}$ is normal in $D$.
 }

 \vspace{3mm}

In 1998, Wang and Fang \cite{w:102} proved

  \vspace{3mm}
 \noindent{\textbf{Theorem D.}}
 \emph{Let $k$, $n\geq k + 1$ be positive integers and $f$ be a transcendental meromorphic function, then $(f^{n})^{(k)}$ assumes every finite
nonzero value infinitely often.
 }

 \vspace{3mm}

 Using the idea of shared values, Li and Gu \cite{L:106} obtained a corresponding normality criteria.

  \vspace{3mm}
 \noindent{\textbf{Theorem E.}}
 \emph{Let $\mathcal {F}$ be a family of meromorphic functions defined in a domain $D$. Let $k$, $n\geq k + 2$ be positive integers and $a\neq0$ be a
finite complex number. If $(f^{n})^{(k)}$ and $(g^{n})^{(k)}$ share a in $D$ for every pair of functions $f$, $g \in\mathcal {F}$, then $\mathcal {F}$ is normal in $D$.
 }
\vspace{3mm}

 In 2004, Alotaibi \cite{A:90} got

\vspace{3mm}
 \noindent{\textbf{Theorem F.}}
 \emph{Suppose that $f$ is a transcendental meromorphic function in the
plane. Let $a\not\equiv0$ be a small function of $f$, then $af(f^{(k)})^{n}-1$ has infinitely many zeros.
 }

 \vspace{3mm}

  Using the idea of shared values, Hu and Meng \cite{H:107} obtained a corresponding normality criteria.

 \vspace{3mm}
 \noindent{\textbf{Theorem G.}}
 \emph{Take positive integers $n$ and $k$ with $n$, $k\geq2$ and take a nonzero complex number $a$. Let $\mathcal{F}$  be a family of meromorphic
functions in the plane domain $D$ such that each $f\in\mathcal {F}$ has only zeros of multiplicity at least $k$. For each pair $(f, g) \in\mathcal {F}$, if $f(f^{(k)})^{n}$ and $g(g^{(k)})^{n}$ share $a$, then $\mathcal {F}$ is normal in $D$.
 }

 \vspace{3mm}
 In 1996, Yang and Hu \cite{Y:98} got

\vspace{3mm}
  \noindent{\textbf{Theorem H.}}
 \emph{Take nonnegative integers $n, n_{1}, \cdots, n_{k}$ with $n\geq 1,\, n_{1}+n_{2}+\cdots+n_{k}\geq1$ and define $d = n +n_{1}+n_{2}+\cdots+n_{k}$.
Let $f$ be a transcendental meromorphic function with the deficiency $\delta(0, f)>3/(3d+1)$. Then for any nonzero value $c$, the function
$f^{n}(f')^{n_{1}}\cdots(f^{(k)})^{n_{k}}-c$ has infinitely many zeros. Moreover, if $n\geq2$, the deficient condition can be omitted.
 }

 \vspace{3mm}
 It's natural to ask whether there exists normality criteria corresponding to Theorem H. We consider this problem and obtain

\begin{Theorem}
Let $a\neq0$ be a constant, $n\geq 2$, $k\geq1$, $n_{k}\geq1$, $n_{j}$ $(j=1, 2, \cdots, k-1)$ be nonnegative integers. Let $\mathcal{F}$  be a family of meromorphic
functions in the plane domain $D$ such that each $f\in\mathcal {F}$ has only zeros of multiplicity at least $k$. For each pair $(f, g) \in\mathcal {F}$, if $f^{n}(f')^{n_{1}}\cdots(f^{(k)})^{n_{k}}$ and $g^{n}(g')^{n_{1}}\cdots(g^{(k)})^{n_{k}}$ share $a$, then $\mathcal {F}$ is normal in $D$.
\end{Theorem}
\begin{Example}
Let $D = \{z: |z| < 1\}$ and $\mathcal {F}= \{f_{m}\}$ where $f_{m} := e^{mz}$, and for every pair of functions $f, g \in\mathcal {F}$, $(f^{n})^{(k)}$ and $(g^{n})^{(k)}$
share $0$ in $D$, it is easy to verify that $\mathcal {F}$ is not normal at the point $z = 0$.
\end{Example}
\begin{Example}
Let $D = \{z: |z| < 1\}$ and $\mathcal {F}= \{f_{m}\}$ where $f_{m} := mz + \frac{1}{mk(k+1)!}$. Then $(f^{k+1}_{m})^{(k)} = m^{k+1}(k + 1)!z + 1$, and for every pair of functions $f, g \in\mathcal {F}$, $(f^{k+1})^{(k)}$ and $(g^{k+1})^{(k)}$
share $1$ in $D$, it is easy to verify that $\mathcal {F}$ is not normal at the point $z = 0$.
\end{Example}
\begin{Remark}
{\rm In Theorem 1.1, let $k=n_{k}=1$, then $f^{n}(f')^{n_{1}}\cdots(f^{(k)})^{n_{k}}=\frac{(f^{n+1})'}{n+1}$ and $g^{n}(g')^{n_{1}}\cdots(g^{(k)})^{n_{k}}=\frac{(g^{n+1})'}{n+1}$, this case is a corollary of Theorem E. Examples 1.1 and 1.2 given by Li and Gu  show that the condition $a\neq 0$ in Theorem E is inevitable and $n\geq k+2$ in Theorem E is sharp. The examples  also show that the conditions in Theorem 1.1 are sharp, at least for the case $k=n_{k}=1$.}
\end{Remark}
If $n=1$, we have
\begin{Theorem}
Let $a\neq0$ be a constant, $k>2$, $n_{k}\geq1$, $n_{j}$ $(j=1, 2, \cdots, k-1)$ be nonnegative integers such that $n_{1}+\cdots+n_{k-1}\geq1$. Let $\mathcal{F}$  be a family of meromorphic functions in the plane domain $D$ such that each $f\in\mathcal {F}$ has only zeros of multiplicity at least $k$. For each pair $(f, g) \in\mathcal {F}$, if $f(f')^{n_{1}}\cdots(f^{(k)})^{n_{k}}$ and $g(g')^{n_{1}}\cdots(g^{(k)})^{n_{k}}$ share $a$, then $\mathcal {F}$ is normal in $D$.
\end{Theorem}
\begin{Remark}
{\rm  In Theorem 1.3, if $n_{1}\geq 2$, the theorem still holds for $k=2$.}
\end{Remark}
 \begin{Corollary}
Let $a\neq0$ be a constant, $n$, $k$, $n_{k}$ be positive integers such that $nk\geq2$ and  $n_{j}$ $(j=1, 2, \cdots, k-1)$ be nonnegative integers. Let $\mathcal{F}$  be a family of holomorphic functions in the plane domain $D$ such that each $f\in\mathcal {F}$ has only zeros of multiplicity at least $k$. For each pair $(f, g) \in\mathcal {F}$, if $f^{n}(f')^{n_{1}}\cdots(f^{(k)})^{n_{k}}$ and $g^{n}(g')^{n_{1}}\cdots(g^{(k)})^{n_{k}}$ share $a$, then $\mathcal {F}$ is normal in $D$.
\end{Corollary}
\begin{Remark}
{\rm  Examples 1.2  shows that Corollary 1.5 fails if $n=k=1$ and thus the condition $nk\geq2$ in Corollary 1.5 is inevitable.}
\end{Remark}
\section{Preliminary lemmas}
\begin{Lemma}[\cite{P:110}] Let $\mathcal {F}$ be a family of meromorphic
functions on the unit disc $\Delta$, all of whose zeros have the
multiplicity at least $k$, and suppose that there exists $A\geq 1$ such
that $|f^{(k)}(z)|\leq A$ wherever $f(z)=0, f\in\mathcal {F}.$ Then
if $\mathcal {F}$ is not normal, there exist, for each $0\leq
\alpha\leq k:$\hfil\break
$(a)$\quad a number $r, \, 0<r<1$,\hfil\break
$(b)$\quad points $z_{n}, |z_{n}|<r,$\hfil\break
$(c)$\quad functions $f_{n}\in \mathcal {F},$ and\hfil\break
$(d)$\quad positive numbers $\rho_{n}\rightarrow 0$\hfil\break
such that
$\rho^{-\alpha}_{n}f_{n}(z_{n}+\rho_{n}\xi)=g_{n}(\xi)\rightarrow
g(\xi)$
 locally uniformly with respect to the spherical metric, where
$g(\xi)$ is a non-constant meromorphic function on $\mathbb{C}$, all of
whose zeros have multiplicity at least $k$, such that $g^{\sharp}(\xi) \leq g^{\sharp}(0)=kA + 1$. In
particular, if $g$ is an entire function, it is of exponential type.
Here, as usual, $g^{\sharp}(z)=|g'(z)|/(1+|g(z)|^{2})$ is the
spherical derivative.
\end{Lemma}
\begin{Lemma}[\cite{C:91}] Let $f$ be an entire function and $M$ a
positive integer. If $f^{\sharp}(z)\leq M$ for all $z\in\mathbb{C},$
then $f$ has the order at most one.
\end{Lemma}
\begin{Lemma} [\cite{Z:91}]
 Take nonnegative integers $n, n_{1}, \cdots, n_{k}$ with $n\geq 1,\,n_{k}\geq1$ and define $d = n +n_{1}+n_{2}+\cdots+n_{k}$. Let $f$ be a transcendental meromorphic function whose zeros have multiplicity at least $k$. Then for any nonzero value $c$, the function
$f^{n}(f')^{n_{1}}\cdots(f^{(k)})^{n_{k}}-c$ has infinitely many zeros, provided that $n_{1}+n_{2}+\cdots+n_{k-1}\geq1$ and $k>2$ if $n=1$. Specially, if $f$ is transcendental entire, the function $f^{n}(f')^{n_{1}}\cdots(f^{(k)})^{n_{k}}-c$ has infinitely many zeros.
\end{Lemma}

\begin{Lemma} Take nonnegative integers $n_{1}, \cdots, n_{k}, k$ with $n_{k}\geq1,\, k\geq2$, $n_{1}+n_{2}+\cdots+n_{k-1}\geq1$ and define $d = 1 +n_{1}+n_{2}+\cdots+n_{k}$.
Let $f$ be a non-constant rational function whose zeros have multiplicity at least $k$. Then for any nonzero value $c$, the function
$f(f')^{n_{1}}\cdots(f^{(k)})^{n_{k}}-c$ has at least two distinct zeros.
\end{Lemma}
{\bf Proof.} We shall divide our argument into two cases.\hfil\break
\noindent{\bf Case 1.} Suppose that $f(f')^{n_{1}}\cdots(f^{(k)})^{n_{k}}-c$ has exactly one zero.\hfil\break
\noindent{\bf Case 1.1.} If $f$ is a non-constant polynomial, since the zeros
of $f$ have multiplicity at least $k$, we know that $f^{n}(f')^{n_{1}}\cdots(f^{(k)})^{n_{k}}$ is also a non-constant polynomial, so $f(f')^{n_{1}}\cdots(f^{(k)})^{n_{k}}-c$ has at least one zero, suppose that
\begin{align}
f(f')^{n_{1}}\cdots(f^{(k)})^{n_{k}}=c+B(z-z_{0})^{l},
\end{align}
where $B$ is a nonzero constant and $l>1$ is an integer. Then (2.5) implies that $f(f')^{n_{1}}\cdots(f^{(k)})^{n_{k}}$ has only simple zeros, which contradicts the assumption that the zeros of $f$ have multiplicity at least $k\geq2$.\hfil\break

\noindent{\bf Case 1.2.} If $f$ is a non-constant rational function but not a polynomial. Set
\begin{align}
f(z)=A \frac{(z-a_{1})^{m_{1}}(z-a_{2})^{m_{2}} \cdots
(z-a_{s})^{m_{s}}}{(z-b_{1})^{l_{1}}(z-b_{2})^{l_{2}} \cdots
(z-b_{t})^{l_{t}}},
\end{align}
where $A$ is a nonzero constant and $m_{i}\geq k\, (i = 1, 2, \cdots
, s), l_{j}\geq 1\, ( j = 1, 2, \cdots , t)$.\hfil\break
Then
\begin{align}
f^{(k)}(z)=A \frac{(z-a_{1})^{m_{1}-k}(z-a_{2})^{m_{2}-k} \cdots
(z-a_{s})^{m_{s}-k}g_{k}(z)}{(z-b_{1})^{l_{1}+k}(z-b_{2})^{l_{2}+k} \cdots
(z-b_{t})^{l_{t}+k}}.
\end{align}
 For simplicity, we denote
\begin{align}
m_{1}+m_{2}+ \cdots +m_{s}=M\geq ks,
\end{align}
\begin{align}
l_{1}+l_{2}+ \cdots +l_{t}=N\geq t.
\end{align}
It is easily obtained that
\begin{align}
\deg(g_{k})\leq k(s+t-1),
\end{align}
moreover, $g_{k}(z)=(M-N)(M-N-1)\cdots(M-N-k+1)z^{k(s+t-1)}+c_{m}z^{k(s+t-1)-1}+\cdots+c_{0}$.
Combining (2.6) and (2.7) yields
\begin{align}
f(f')^{n_{1}}\cdots(f^{(k)})^{n_{k}}=A^{d} \frac{(z-a_{1})^{dm_{1}-\sum_{j=1}^{k}jn_{j}}
\cdots(z-a_{s})^{dm_{s}-\sum_{j=1}^{k}jn_{j}}g(z)}{(z-b_{1})^{dl_{1}+\sum_{j=1}^{k}jn_{j}}
\cdots (z-b_{t})^{dl_{t}+\sum_{j=1}^{k}jn_{j}}}=\frac{P(z)}{Q(z)},
\end{align}
where $P(z)$, $Q(z)$, $g(z)$ are polynomials and  $g(z)=\prod_{j=1}^{k}g_{j}^{n_{j}}(z)$ with $\deg(g)\leq\sum_{j=1}^{k}jn_{j}(s+t-1)$.
Moreover $g(z)=(M-N)^{d-1}(M-N-1)^{d-1-n_{1}}\cdots(M-N-k+1)^{n_{k}}z^{\sum_{j=1}^{k}jn_{j}(s+t-1)}+d_{m}z^{\sum_{j=1}^{k}jn_{j}(s+t-1)-1}+\cdots+d_{0}$. \hfil\break
Then from (2.10) we obtain
\begin{align}
(f(f')^{n_{1}}\cdots(f^{(k)})^{n_{k}})'=A^{d} \frac{(z-a_{1})^{dm_{1}-\sum_{j=1}^{k}jn_{j}-1}
\cdots(z-a_{s})^{dm_{s}-\sum_{j=1}^{k}jn_{j}-1}h(z)}{(z-b_{1})^{dl_{1}+\sum_{j=1}^{k}jn_{j}+1}
\cdots (z-b_{t})^{dl_{t}+\sum_{j=1}^{k}jn_{j}+1}},
\end{align}
where $h(z)$ is a polynomial with $\deg(h)\leq (\sum_{j=1}^{k}jn_{j}+1)(s+t-1)$. \vskip 2pt Since
$f(f')^{n_{1}}\cdots(f^{(k)})^{n_{k}}-c$ has exactly one zero, we obtain from (2.11) that
\begin{align}
f(f')^{n_{1}}\cdots(f^{(k)})^{n_{k}}=c+\frac{B(z-z_{0})^{l}}{(z-b_{1})^{dl_{1}+\sum_{j=1}^{k}jn_{j}}
\cdots (z-b_{t})^{dl_{t}+\sum_{j=1}^{k}jn_{j}}},
\end{align}
where $B$ is a nonzero constant.
Then
\begin{align}
(f(f')^{n_{1}}\cdots(f^{(k)})^{n_{k}})'=\frac{B(z-z_{0})^{l-1}H(z)}{(z-b_{1})^{dl_{1}+\sum_{j=1}^{k}jn_{j}+1}
\cdots (z-b_{t})^{dl_{t}+\sum_{j=1}^{k}jn_{j}+1}},
\end{align}
where $H(z)$ is a polynomial of the form
\begin{align}
H(z)=B(l-dN-\sum_{j=1}^{k}jn_{j}t)z^{t}+B_{t-1}z^{t-1}+\cdots+b_{0}.
\end{align}
\noindent{\bf Case 1.2.1.} If $l\neq dN+\sum_{j=1}^{k}jn_{j}t$.\hfil\break
From (2.11) we get $\deg(p)\geq \deg(Q)$, namely
\begin{align}
dM-\sum_{j=1}^{k}jn_{j}s+\deg(g)\geq dN+\sum_{j=1}^{k}jn_{j}t.
\end{align}
In view of $\deg(g)\leq\sum_{j=1}^{k}jn_{j}(s+t-1)$, we deduce from (2.16) that
\begin{align}
d(M-N)\geq \sum_{j=1}^{k}jn_{j},
\end{align}
which implies $M>N$.\hfil\break
Combining (2.12), (2.14) and (2.15) yields
\begin{align}
dM-(\sum_{j=1}^{k}jn_{j}+1)s\leq\deg(H)=t\leq N< M,
\end{align}
which implies that
\begin{align}
\sum_{j=1}^{k}n_{j}M-\sum_{j=1}^{k}jn_{j}s< s,
\end{align}
this together with (2.8) implies that
\begin{align}
 \sum_{j=1}^{k-1}(k-j)n_{j}<1.
\end{align}
which is a contradiction since $n_{1}+n_{2}+\cdots+n_{k-1}\geq1$.\hfil\break

\noindent{\bf Case 1.2.2.} If $l= dN+\sum_{j=1}^{k}jn_{j}t$.\hfil\break
If $M>N$, with similar discussion as above, we get the same contradiction.\hfil\break
If $M\leq N$, combining (2.12) and (2.14) yields
\begin{align}
 l-1\leq \deg(h) \leq (\sum_{j=1}^{k}jn_{j}+1)(s+t-1),
\end{align}
we deduce from (2.21) that
\begin{align}
 dN \leq \sum_{j=1}^{k}jn_{j}s+s+t-\sum_{j=1}^{k}jn_{j}< \sum_{j=1}^{k}jn_{j}s+s+N.
\end{align}
Note that $M\leq N$, from (2.22) we get
\begin{align}
 (d-1)M -\sum_{j=1}^{k}jn_{j}s<s.
\end{align}
Similar to the end of the proof in Case 1.2.1, we get a contradiction. Case 1 has been ruled out.\hfil\break

\noindent{\bf Case 2.} Suppose that $f(f')^{n_{1}}\cdots(f^{(k)})^{n_{k}}-c$ has no zero.\hfil\break
Similar to the proof of Case 1.1, we deduce that $f$ is not a polynomial. If $f$ is a non-constant rational function but not a polynomial.
Similar to the proceeding of proof in Case 1.2.1, we get a contradiction.\hfil\break
Hence $f^{n}(f')^{n_{1}}\cdots(f^{(k)})^{n_{k}}-c$ has at least two distinct zeros.\hfil\break
 This proves Lemma 2.4.

 \vspace{3mm}
 Using the similar proof of Lemma 2.4, we get
\begin{Lemma}
 Take nonnegative integers $n$, $n_{1}, \cdots, n_{k}, k$ with $n\geq2$, $n_{k}\geq1,\, k\geq1$.
Let $f$ be a non-constant rational function whose zeros have multiplicity at least $k$. Then for any nonzero value $c$, the function
$f^{n}(f')^{n_{1}}\cdots(f^{(k)})^{n_{k}}-c$ has at least two distinct zeros.
\end{Lemma}
\begin{Lemma}
 Take nonnegative integers $n$, $n_{1}, \cdots, n_{k}, k$ with $n\geq1$, $n_{k}\geq1,\, k\geq2$.
Let $f$ be a non-constant polynomial whose zeros have multiplicity at least $k$. Then for any nonzero value $c$, the function
$f^{n}(f')^{n_{1}}\cdots(f^{(k)})^{n_{k}}-c$ has at least two distinct zeros.
\end{Lemma}

\section{Proof of Theorem 1.3}
Without loss of generality, we may assume $D = \Delta = \{z: |z| < 1\}$. For simplicity, we denote $f_{j}(z_{j}+\rho_{j}\xi)$ by $f_{j}$ and set
\begin{align*}
\gamma_{M}=1+n_{1}+\cdots+n_{k}, \quad \Gamma_{M}=\sum_{j=1}^{k}jn_{j}
\end{align*}
Suppose that $\mathcal {F}$ is not normal in $D$. By Lemma 2.1, for $0\leq \alpha <k$, there
exist:  $r<1$, $z_{j}\rightarrow 0$ $(j\rightarrow\infty)$, $f_{j}\in \mathcal {F}$ and $\rho_{j}\rightarrow 0^{+}$
 such that $g_{j}(\xi)=\rho^{-\Gamma_{M}/\gamma_{M}}_{j}f_{j}(z_{j}+\rho_{j}\xi)\rightarrow
g(\xi)$ locally uniformly with respect to the spherical metric, where
$g(\xi)$ is a non-constant meromorphic function on $\mathbb{C}$, all of
whose zeros have multiplicity at least $k$, such that $g^{\sharp}(\xi) \leq g^{\sharp}(0)=kA + 1$.\hfil\break
On every compact subset of $\mathbb{C}$ which contains no poles of $g$, we have uniformly
\begin{eqnarray}
f_{j}(f_{j}')^{n_{1}}\cdots(f_{j}^{(k)})^{n_{k}}-a&=&g_{j}(\xi)(g_{j}'(\xi))^{n_{1}}\cdots(g_{j}^{(k)}(\xi))^{n_{k}}-a\nonumber\\&\rightarrow &g(g')^{n_{1}}\cdots(g^{(k)})^{n_{k}}-a.
\end{eqnarray}
If $g^{n}(g')^{n_{1}}\cdots(g^{(k)})^{n_{k}}\equiv a$, then $g$ has no zeros and no poles, thus $g$ is an entire function. By Lemma 2.2, $g$ is of order at most 1. Moreover, $g(\xi)=e^{c_{1}\xi+c_{0}}$, where $c_{1}(\neq 0)$ and $c_{0}$ are constants. Thereby, we get
\begin{align}
g(\xi)(g'(\xi))^{n_{1}}\cdots(g^{(k)}(\xi))^{n_{k}}=c_{1}^{\Gamma_{M}}e^{\gamma_{M}(c_{1}\xi+c_{0})},
\end{align}
which contradicts the case $g^{n}(g')^{n_{1}}\cdots(g^{(k)})^{n_{k}}\equiv a$.\hfil\break
Since $g$ is a non-constant meromorphic function, by Lemmas 2.3 and 2.4, we deduce that $g(\xi)(g'(\xi))^{n_{1}}\cdots(g^{(k)}(\xi))^{n_{k}}-a$ has at least
two distinct zeros.\hfil\break
We claim that $g(\xi)(g'(\xi))^{n_{1}}\cdots(g^{(k)}(\xi))^{n_{k}}-a$ has just a unique zero.\hfil\break
Let $\xi_{0}$ and $\xi_{0}^{\star}$ be two distinct zeros of $g(\xi)(g'(\xi))^{n_{1}}\cdots(g^{(k)}(\xi))^{n_{k}}-a$. We choose a positive number $\delta$ small enough such that $D_{1} \bigcap D_{2} =\emptyset$, and such that $g(g')^{n_{1}}\cdots(g^{(k)})^{n_{k}}-a$ has no other zeros in $D_{1}\bigcup D_{2}$
except for $\xi_{0}$ and $\xi_{0}^{\star}$, where
\begin{align*}
D_{1}=\{\xi \in \mathbb{C}\mid |\xi-\xi_{0}|<\delta \}, \quad D_{2}=\{\xi \in \mathbb{C}\mid |\xi-\xi_{0}^{\star}|<\delta \}.
\end{align*}
By (3.1) and Hurwitz's theorem, for sufficiently large $j$ there exist points $\xi_{j}\in D_{1}$, $\xi_{j}^{\star}\in D_{2}$ such that
\begin{align*}
f_{j}(z_{j}+\rho_{j}\xi_{j})(f'_{j}(z_{j}+\rho_{j}\xi_{j}))^{n_{1}}\cdots(f_{j}^{(k)}(z_{j}+\rho_{j}\xi_{j}))^{n_{k}}-a=0,
\end{align*}
\begin{align*}
f_{j}(z_{j}+\rho_{j}\xi_{j}^{\star})(f'_{j}(z_{j}+\rho_{j}\xi_{j}^{\star}))^{n_{1}}\cdots(f_{j}^{(k)}(z_{j}+\rho_{j}\xi_{j}^{\star}))^{n_{k}}-a=0.
\end{align*}
By the assumption in Theorem 1.1, $f_{1}(f'_{1})^{n_{1}}\cdots(f_{1}^{(k)})^{n_{k}}$ and $f_{j}(f'_{j})^{n_{1}}\cdots(f_{j}^{(k)})^{n_{k}}$ share $a$ for each $j$, it follows that
\begin{align*}
f_{1}(z_{j}+\rho_{j}\xi_{j})(f'_{1}(z_{j}+\rho_{j}\xi_{j}))^{n_{1}}\cdots(f_{1}^{(k)}(z_{j}+\rho_{j}\xi_{j}))^{n_{k}}-a=0,
\end{align*}
\begin{align*}
f_{1}(z_{j}+\rho_{j}\xi_{j}^{\star})(f'_{1}(z_{j}+\rho_{j}\xi_{j}^{\star}))^{n_{1}}\cdots(f_{1}^{(k)}(z_{j}+\rho_{j}\xi_{j}^{\star}))^{n_{k}}-a=0.
\end{align*}
Let $j\rightarrow\infty$, and note that $z_{j}+\rho_{j}\xi_{j}\rightarrow 0$, $z_{j}+\rho_{j}\xi_{j}^{\star}\rightarrow 0$, we get
\begin{align*}
f_{1}(0)(f'_{1}(0))^{n_{1}}\cdots(f_{1}^{(k)}(0))^{n_{k}}-a=0.
\end{align*}
Since the zeros of $f_{1}(f'_{1})^{n_{1}}\cdots(f_{1}^{(k)})^{n_{k}}-a$ have no accumulation points, in fact, for sufficiently large $j$, we have
\begin{align*}
z_{j}+\rho_{j}\xi_{j}=0, \quad z_{j}+\rho_{j}\xi_{j}^{\star}=0.
\end{align*}
Thus
\begin{align*}
\xi_{j}=-\frac{z_{j}}{\rho_{j}}, \quad \xi_{j}^{\star}=-\frac{z_{j}}{\rho_{j}}.
\end{align*}
This contradicts the fact that\hfil\break
 $\xi_{j}\in D(\xi_{0}, \delta)$, $\xi_{j}^{\star}\in D(\xi_{0}^{\star}, \delta)$ and $D(\xi_{0}, \delta)\bigcap D(\xi_{0}^{\star}, \delta)=\emptyset$.\hfil\break
  So $g(\xi)(g'(\xi))^{n_{1}}\cdots(g^{(k)}(\xi))^{n_{k}}-a$ has just a unique zero, which contradicts the fact that $g(\xi)(g'(\xi))^{n_{1}}\cdots(g^{(k)}(\xi))^{n_{k}}-a$ has at least
two distinct zeros.\hfil\break
This completes the proof of Theorem 1.3.

\vspace{3mm}
By Theorem H, Lemmas 2.3, 2.5 and 2.6, the proofs of Theorem 1.1 and Corollary 1.5 can be carried out in the line of Theorem 1.3, we omit the
process here.
\section{Discussion}
In Theorem 1.3, if $n_{1}=n_{2}= \cdots =n_{k-1}=0$, then $f(f')^{n_{1}}\cdots(f^{(k)})^{n_{k}}=f(f^{(k)})^{n_{k}}$, $g(g')^{n_{1}}\cdots(g^{(k)})^{n_{k}}=g(g^{(k)})^{n_{k}}$. This case is the same as the case in Theorem G. Hu and Meng proved that Theorem G holds for the case $k\geq2, n_{k}\geq2$. They also gave an example \cite[Example 1.3]{H:107}  to show that Theorem G is not valid for the case $k=1$. It's natural to ask whether Theorem 1.3 holds for the case $n_{k}=1$, $n_{1}+\cdots+n_{k-1}=0$. Actually, this is an open problem as follows.
\begin{Problem}
Let $\mathcal {F}$ be a family of meromorphic functions in a domain $D$, let $k$ be a positive integer and $b$ be a finite nonzero value.
 If, for every $f\in \mathcal {F}$, all zeros of $f$ have multiplicity at least $k$, and $f(z)f^{(k)}(z)\neq b$, is $\mathcal {F}$
normal in $D$?
\end{Problem}
Xu and Cao \cite{X:106} gave a partial answer to Problem 4.1.
\begin{Theorem}
Let $\mathcal {F}$ be a family of meromorphic functions in a domain $D$, let $k$ be a positive integer and $b$ be a finite nonzero value.
 If, for every $f\in \mathcal {F}$, all zeros of $f$ have multiplicity at least $k+1$, and $f(z)f^{(k)}(z)\neq b$, then $\mathcal {F}$ is
normal in $D$.
\end{Theorem}

\end{document}